\documentclass[12pt]{article}
\usepackage[cp866]{inputenc}
\usepackage{amsmath}
\usepackage{amsfonts}
\usepackage{amssymb}
\usepackage{longtable}

\newtheorem{theorem}{Theorem}

\newcommand{\spa}{\mathrm{span}}
\newcommand{\pr}{\mathrm{pr}}
\newcommand{\Real}{\mathbb{R}}
\newcommand{\g}{\mathfrak{g}}
\newcommand{\h}{\mathfrak{h}}
\newcommand{\so}{\mathfrak{so}}

\newcommand{\gl}{\mathfrak{gl}}
\newcommand{\zr}{\ltimes}
\newcommand{\hol}{\mathfrak{hol}}
\newcommand{\Ga}{\Gamma}

\begin{document}



\author{A.\,S.~Galaev}

\title{Note on the holonomy groups of pseudo-Riemannian manifolds}


\maketitle

\begin{abstract}
For an arbitrary
subalgebra~$\mathfrak{h}\subset\mathfrak{so}(r,s)$, a polynomial
pseudo-Riemannian metric of signature~$(r+2,s+2)$ is constructed,
the holonomy algebra of this metric contains~$\mathfrak{h}$ as a
subalgebra. This result shows the essential distinction of the
holonomy algebras of pseudo-Riemannian manifolds of index bigger
or equal to~2 from the holonomy algebras of Riemannian and
Lorentzian manifolds.

\end{abstract}


\section{Introduction} \label{subsec1}

The holonomy group of a linear connection on an~$n$-dimensional
manifold is contained in the Lie group~${\rm GL}(n,\Real)$ and it
represents an important invariant of the connection. Hano and
Ozeki~\cite{H-O} showed that any connected linear Lie
group~$G\subset{\rm GL}(n,\Real)$ may be realized as the holonomy
group of a space with a linear connection. This connection,  as a
rule, are of non-zero torsion. On the contrary, the absence of the
torsion imposes some algebraic condition on the holonomy group.
Berger used this condition  in order to get the classification of
the connected holonomy groups of Riemannian manifolds and of
connected irreducible holonomy groups of the spaces with
torsion-free linear connections~\cite{Ber}. Later the Berger lists
were corrected. The classification of the connected irreducible
holonomy groups of torsion-free linear connections is obtained
in~\cite{M-Sch}. The most important result turned out to be the
classification of the connected holonomy groups of Riemannian
manifolds, this result has many applications in geometry and
theoretical physics~\cite{Al,Besse,Joyce07,Ego11}. Lately the
interest to the pseudo-Riemannian manifolds appears. Recently the
classification of the connected holonomy groups of Lorentzian
manifolds was obtained~\cite{ESI,Leistner}. Note that the study of
the connected holonomy groups is equivalent to the study of the
holonomy algebras, i.e. the corresponding Lie algebras. Consider a
Lorentzian manifold of dimension~$n+2\geq 4$. Its holonomy algebra
is contained in the Lorentzian Lie algebra~$\so(1,n+1)$. It is
enough to consider the holonomy algebras contained in the maximal
algebra preserving an isotropic line in the Minkowski
space,~$\g\subset\so(1,n+1)_\ell=(\Real\oplus\so(n))\zr\Real^n$.
With the algebra~$\g$ its projection~$\h$ to~$\so(n)$ is
associated. The key fact is that~$\h\subset\so(n)$ is the holonomy
algebra of a Riemannian manifold.

Note that the maximal subalgebra of the pseudo-orthogonal Lie
algebra~$\so(r+2,s+2)$, $r+s\geq 2$, preserving a two-dimensional
isotropic subspace of the pseudo-Euclidean space~$\Real^{r+2,s+2}$
has the
form~$(\gl(2,\Real)\oplus\so(r,s))\zr(\Real^2\otimes\Real^{r,s}\zr\Real)$.
In this note, for any
subalgebra~$\mathfrak{h}\subset\mathfrak{so}(r,s)$ a polynomial
pseudo-Riemannian metric of signature~$(r+2,s+2)$  is constructed,
the holonomy algebra of this metric
equals~$\h\zr(\Real^2\otimes\Real^{r,s}\zr\Real)$. It turns out
that the  holonomy algebra of a pseudo-Riemannian manifold of
signature~$(r+2,s+2)$, i.e. of index bigger or equal to~2, can
depend on an arbitrary
subalgebra~$\mathfrak{h}\subset\mathfrak{so}(r,s)$. This indicates
the non-visibility of the holonomy algebras of pseudo-Riemannian
manifolds of index bigger or equal to~2. In particular, the
holonomy algebra of a pseudo-Riemannian manifold of
signature~$(2,n+2)$, $n\geq 2$, may depend on an arbitrary
subalgebra~$\h\subset\so(n)$; this shows the fundamental
difference from the case of Lorentzian manifolds. In some sense,
the obtained result is analogous to the result from~\cite{H-O}.
Other results on the holonomy groups of pseudo-Riemannian
manifolds can be found in the review~\cite{IRMA}.

\section{Results} \label{subsec2}

Consider the pseudo-Euclidean space~$\Real^{r+2,s+2}$ with the
pseudo-Euclidean metric~$\eta$ of signature~$(r+2,s+2)$,
$r+s=n\geq 2$. Let us fix a basis~ $p_1,p_2,e_1,...,e_n,q_1,q_2$
of the space~$\Real^{r+2,s+2}$ such that~$\eta$ has the following
non-zero values:
$$\eta(p_1,q_1)=\eta(p_2,q_2)=1,\quad
\eta(e_i,e_i)=\epsilon_i,\quad
\epsilon_1=\cdots=\epsilon_r=-1,\quad
\epsilon_{r+1}=\cdots=\epsilon_n=1.$$ Consider the
subalgebra~$\so(r+2,s+2)_{<p_1,p_2>}\subset\so(r+2,s+2)$
preserving the isotropic
subspace~$\spa\{p_1,p_2\}\subset\Real^{r+2,s+2}$. This Lie algebra
has the following matrix form:
  $$\left\{\left. \left (\begin{array}{ccc}
B & \begin{array}{c}  -(E_{r,s}X)^t \\ -(E_{r,s}Y)^t \end{array} &\begin{array}{cc}  0&-c \\ c&0 \end{array}  \\
\begin{array}{cc} 0&0\end{array} &  A& \begin{array}{cc} X&Y\end{array}\\
\begin{array}{cc}  0&0 \\
0&0 \end{array} & \begin{array}{c}  0 \\ 0 \end{array} & -B^t
\end{array}\right)\right|\,\begin{array}{ll} B\in\gl(2,\Real),& A\in\so(r,s),\\X,Y\in\Real^{r,s},&c\in \Real\end{array}\right\},$$
where $E_{r,s}={\rm diag}(\epsilon_1,...,\epsilon_n)$. We get the
decomposition
$$\so(r+2,s+2)_{<p_1,p_2>}=(\gl(2,\Real)\oplus\so(r,s))\zr(\Real^2\otimes\Real^{r,s}\zr\Real),$$
where the symbol~$\zr$ indicates that on the right side from it is
an ideal.

For an arbitrary subalgebra~$\h\subset\so(r,s)$ we define the
subalgebra
$$\g^\h=\left\{\left.\left (\begin{array}{ccccc}
0&0 &-(E_{r,s}X)^t &0 &-c\\
0&0&-(E_{r,s}Y)^t &c &0\\
0 &0&A&X&Y\\
0&0&0&0&0\\
0&0&0&0&0\\
\end{array}\right)\right|\, A\in\h,X,Y\in\Real^{r,s},c\in\Real\right\}$$ of the Lie algebra~$\so(2,n+2)_{<p_1,p_2>}.$
We denote an element of the Lie algebra~$\g^\h$ by~$(A,X,Y,c)$.
The non-zero Lie brackets are the following:
\begin{align}\label{LieB1} [(A,0,0,0),(A_1,X,Y,0)]&=([A,A_1],AX,AY,0),\\ \label{LieB2}
[(0,X,0,0),(0,0,Y,0)]&=(0,0,0,\eta(X,Y)).\end{align} We get the
decomposition
$$\g^\h=\h\zr(\Real^2\otimes\Real^{r,s}\zr\Real).$$

\begin{theorem} For any subalgebra~$\h\subset\so(r,s)$, the algebra~$\g^\h$ is the holonomy algebra of a
pseudo-Riemannian manifold of signature~$(r+2,s+2)$.
\end{theorem}

{\bf Proof.} For an arbitrary subalgebra~$\h\subset\so(r,s)$ we
construct a pseudo-Riemannian metric on the space~$\Real^{r+s+4}$
and show that the holonomy algebra of this metric at the point~$0$
coincides with~$\g^\h$.

{\bf Construction of the metric.} We fix  elements
~$B_1,...,B_N\in \h$, generating the Lie algebra~$\h$. Consider
the matrices~$(B^i_{j\alpha})^n_{i,j=1}$ of these elements with
respect to the basis~$e_1,...,e_n$ of the space~$\Real^{r,s}$. The
condition~$B_\alpha\in\h\subset\so(r,s)$ implies
$$\epsilon_jB^i_{j\alpha}=-\epsilon_iB^j_{i\alpha}.$$

Let~$v,z,x^1,...,x^n,u,w$ be the standard coordinates
on~$M=\Real^{r+s+4}$. Consider the pseudo-Riemannian metric
$$g=2dvdu+2dzdw+\sum_{i=1}^n\epsilon_i(dx^i)^2+\sum_{i=1}^n2A_idx^idw+f(du)^2+f(dw)^2,$$
where
$$A_i=\sum_{\alpha=1}^N\sum_{j=1}^nA_{ij\alpha}x^ju^\alpha,\quad
A_{ij\alpha}=\epsilon_iB^i_{j\alpha}, \quad
f=\sum_{i=1}^n(x^i)^2.$$ The metric~$g$ is of
signature~$(r+2,s+2)$. Let~$\hol$ be the holonomy  algebra of this
metric at the point~0. Consider the basis
$$p_1=(\partial_v)_0,\,p_2=(\partial_z)_0,\,e_1=(\partial_1)_0,\dots,
e_n=(\partial_n)_0,\,q_1=(\partial_u)_0,\,q_2=(\partial_w)_0$$ of
the tangent space~$T_0M$. We get~$\eta=g_0$, this allows us to
identify$(T_0M,g_0)$ with~$(\Real^{r+2,s+2},\eta)$ and the
holonomy algebra~$\hol$ with a subalgebra of~$\so(r+2,s+2)$.

{\bf Computation of the holonomy algebra.} The constructed metric
is analytical. From the proof of theorem~9.2 from~\cite{K-N1} it
follows that~$\hol$ is generated by the elements of the form
$$\nabla_{\partial_{a_\alpha}}\cdots\nabla_{\partial_{a_1}}R(\partial_a,\partial_b)(0)\in\so(T_0M,g_0)=\so(r+2,s+2),
\quad \alpha=0,1,2...,$$
 where~$\nabla$ is the Levi-Civita connection defined by the
 metric~ $g$ and~$R$ is the curvature tensor. The indices~$a,b,c$ will run through all coordinates on~$M$,
the indices~$i,j,k$  will take the values~$1,...,n$.
  The Levi-Civita connection is determined by  its Christoffel symbols~$\Gamma^a_{bc}$,
 $\nabla_{\partial_b}\partial_c=\sum_a\Gamma^a_{bc}\partial_a,$
 which can be found using the formula $$\Gamma^a_{bc}=\Gamma^a_{cb}=\frac{1}{2}
 \sum_d g^{ad}(\partial_cg_{bd}+\partial_bg_{cd}-\partial_dg_{bc}),$$
 where~$(g^{ab})$ is the inverse matrix to the matrix~$(g_{ab})$ of the metric~$g$. The components of the
 curvature tensor are defined by the equality~$R(\partial_a,\partial_b)\partial_c=\sum_dR^d_{cab}\partial_d$
and can be found in the following way:
$$R^d_{cab}=\partial_a\Gamma^d_{bc}-\partial_b\Gamma^d_{ac}+
  \sum_e(\Ga^e_{bc}\Ga^d_{ae}-\Ga^e_{ac}\Ga^d_{be}).$$
 For the matrix~$(g^{ab})$ we get
 $$(g^{ab})=\left(\begin{array}{ccc}-fE_2&C&E_2\\C^t&E_{r,s}&0\\E_2&0&0\end{array}\right),\quad
 C=\left(\begin{array}{ccc}0&\cdots&0\\-\epsilon_1A_1&\cdots&-\epsilon_nA_n\end{array}\right).$$
 In  order to find the holonomy algebra we will need only the following Christoffel symbols:
 \begin{align}
 \Ga^i_{j k}&=0, \quad  \Ga^i_{j u}=0, \quad
 \Ga^i_{j w}=\sum_{\alpha=1}^NB^i_{j \alpha}u^\alpha, \label{e21}\\
  \Ga^{a}_{v b}&=\Ga^{a}_{z b}=0,\quad \Ga^u_{ia}=\Ga^w_{ia}=0, \label{e20}\\
 \Ga^{v}_{ij}&=\Ga^{z}_{ij}=0,\quad \Ga^v_{iu}=x^i,\quad
 \Ga^z_{iw}=x^i+\sum_k\epsilon_kA_k(\partial_iA_k-\partial_kA_i),  \label{e25}
\end{align}
 the following components of the curvature tensor:
\begin{equation}
R^{i}_{j uw}=\sum_{\alpha=1}^N\alpha B^i_{j
\alpha}u^{\alpha-1},\quad R^{i}_{j a b}=0, \text{ if }
\{a\}\cup\{b\}\neq\{u,w\}, \label{e50}\end{equation} and  the
following components of the curvature tensor at the point~$0$:
\begin{equation} R^{v}_{i i u}(0)=1, \quad R^{z}_{i i w}(0)=1.
\label{e60}\end{equation} The computation of these values are
direct. Note that there exists the following recurrent formula:
\begin{equation}\label{nabnabR}\nabla_{a_\alpha}\cdots\nabla_{a_1}R^d_{cab}=
\partial_{a_\alpha}\nabla_{a_{\alpha-1}}\cdots\nabla_{a_1}R^d_{cab}+
[\Ga_{a_\alpha},\nabla_{\partial_{a_{\alpha-1}}}\cdots\nabla_{\partial_{a_1}}R(\partial_a,\partial_b)]^d_c,
\end{equation} where~$\Gamma_{a_\alpha}$ denote the operator with the matrix~$(\Ga_{ba_\alpha}^a)$.

{\bf Proof of the inclusion~$\hol\subset\g^\h$.} The
equality~$\Ga^{a}_{v b}=\Ga^{a}_{z b}=0$ means that the vector
fields~$\partial_v,\partial_z$ are parallel,
i.e.~$\nabla\partial_v=\nabla\partial_z=0$. According to the
holonomy principle, $\hol$~annihilates the
vectors~$p_1=(\partial_v)_0$,~$p_2=(\partial_z)_0$, this implies
the inclusion~$\hol\subset\g^{\so(r,s)}$. It remains to prove
that~$\pr_{\so(r,s)}\hol\subset\h$, i.e.
$$\pr_{\so(r,s)}(\nabla_{\partial_{a_\alpha}}\cdots\nabla_{\partial_{a_1}}R(\partial_a,\partial_b)(0))\in
\h$$ for all~$\alpha$. The equalities~\eqref{LieB2}
and~\eqref{e21} show that
\begin{multline*}\nabla_{a_\alpha}\cdots\nabla_{a_1}R^i_{jab}=
\partial_{a_\alpha}\nabla_{a_{\alpha-1}}\cdots\nabla_{a_1}R^i_{jab}\\+
[\pr_{\so(r,s)}\Gamma_{a_\alpha},\pr_{\so(r,s)}(\nabla_{\partial_{a_{\alpha-1}}}
\cdots\nabla_{\partial_{a_1}}R(\partial_a,\partial_b))]^i_j.\end{multline*}
Note that if~$\pr_{\so(r,s)}\Gamma_{a_\alpha}\neq 0$,
then~$a_\alpha=w$; in this case
$$\pr_{\so(r,s)}\Gamma_w=\sum_{\alpha=1}^NB_\alpha u^\alpha.$$
Now it is easy to prove the inclusion~$\hol\subset\g^{\h}$ using
the induction and~\eqref{e21},~\eqref{e50}.

{\bf Proof of the inclusion~$\g^\h\subset\hol$}.
Equalities~\eqref{e60} imply the inclusion\\
$\{(0,X,Y,0)|X,Y\in\Real^{r,s}\}\subset\hol$. Using~\eqref{LieB2},
we get~$\{(0,0,0,c)|c\in\Real\}\subset\hol$. From~\eqref{e21}
and~\eqref{nabnabR} follows the equality
$$(\nabla_u)^{\beta}R^i_{juw}=(\partial_u)^{\beta}R^i_{juw}=\sum_{\alpha=\beta+1}^N
\alpha(\alpha-1)\cdots(\alpha-\beta)B^i_{j\alpha}u^{\alpha-\beta-1}$$
for all~$\quad 0\leq \beta\leq N-1.$ In particular,
$$(\nabla_u)^{\beta}R^i_{juw}(0)=(\beta+1)!B^i_{j\beta+1},\quad 0\leq \beta\leq N-1,$$
i.e.
$$\pr_{\so(r,s)}(\nabla_{\partial_u})^{\beta-1}R(\partial_u,\partial_w)(0)=\beta!B_{\beta},\quad
1\leq \beta\leq N.$$ Since the Lie algebra~$\h$ is generated by
the elements~$B_1,...,B_N$, the last equality implies the
inclusion~$\g^\h\subset\hol$. The theorem is true.

\section{Correlation with the case of Lorentzian manifolds} Let
us compare the obtained result with the classification of the
holonomy algebras of Lorentzian manifolds \cite{ESI,Leistner}. The
holonomy algebra of an~$(n+2)$-dimensional Lorentzian manifold is
a subalgebra of the Lorentzian Lie algebra~$\so(1,n+1)$, $n\geq
0$. Consider a basis~$p,e_1,...,e_n,q$ of the Minkowski
space~$\Real^{1,n+1}$ such that the Minkowski metric~$\eta$ has
only the following non-zero values: $\eta(p,q)=\eta(e_i,e_i)=1$.
The subalgebra~$\so(1,n+1)_{\Real p}\subset\so(1,n+1)$ preserving
the isotropic line~$\Real p$ has the form
$$\left\{\left. \left (\begin{array}{ccc} a
&-X^t & 0\\ 0 & A &X \\ 0 & 0 & -a \\
\end{array}\right)\right|\, a\in \Real,\, \,A \in \so(n),\, X\in \Real^n
\right\}=(\Real\oplus\so(n))\zr\Real^n.$$ It is enough to consider
the subalgebras~$\g\subset\so(1,n+1)$ contained
in~$\so(1,n+1)_{\Real p}$. For an arbitrary
subalgebra~$\h\subset\so(n)$ consider the subalgebra
$$\g^\h=\left\{\left. \left (\begin{array}{ccc} 0
&-X^t & 0\\ 0 & A &X \\ 0 & 0 & 0 \\
\end{array}\right)\right|\, A \in \h,\, X\in \Real^n
\right\}=\h\zr\Real^n\subset\so(1,n+1)_{\Real p}.$$
Leistner~\cite{Leistner} proved the following non-trivial
statement: {\it if~$\g^\h$ is the holonomy algebra of a Lorentzian
manifold, then the subalgebra~$\h\subset\so(n)$ must be the
holonomy algebra of a Riemannian manifold.} The proof is based on
the fact that the holonomy algebra~$\g\subset\so(r,s)$ of an
arbitrary pseudo-Riemannian manifold of signature~$(r,s)$ is
generated by the images of algebraic curvature tensors, these
tensors belong to the space~$\mathcal{R}(\g)$ consisting of the
2-forms on~$\Real^{r,s}$ with the values in~$\g$ and satisfying
the first Bianchi identity
$$R(X,Y)Z+R(Y,Z)X+R(Z,X)Y=0,\quad X,Y,Z\in\Real^{r,s}.$$ For~$R\in\mathcal{R}(\g^\h)$,
the projection~$\pr_\h\circ R$ is of the form
$$\pr_\h\circ R(p,\cdot)=0,\quad \pr_\h\circ
R|_{\Real^n\times\Real^n}\in\mathcal{R}(\h),\quad \pr_\h\circ
R(X,q)=P(X),\quad X\in\Real^n,$$ where $P:\Real^n\to\h$ is a
linear map satisfying the identity
$$\eta(P(X)Y,Z)+\eta(P(Y)Z,X)+\eta(P(Z)X,Y)=0,\quad X,Y,Z\in\Real^n.$$
Let~$\mathcal{P}(\h)$ be the space of such maps~$P$. We get
that~$\h$ must be generated by the images of the elements of the
spaces~$\mathcal{R}(\h)$ and~$\mathcal{P}(\h)$. Leistner showed
that from this condition it follows that~$\h$ must be generated by
the elements of the space~$\mathcal{R}(\h)$; this means
that~$\h\subset\so(n)$ is the holonomy algebra of a Riemannian
manifold.

Consider now the subalgebra~$\g^\h\subset\so(2,n+2)$ from the
previous section, where~$\h\subset\so(n)$. It is easy to show that
for~$R\in\mathcal{R}(\g^\h)$, the projection~$\pr_\h\circ R$
satisfies
$$\pr_\h\circ R(p_1,\cdot)=\pr_\h\circ R(p_1,\cdot)=0,\quad \pr_\h\circ
R|_{\Real^n\times\Real^n}\in\mathcal{R}(\h),$$ $$\pr_\h\circ
R(\cdot|_{\Real^n},q_1),\,\pr_\h\circ
R(\cdot|_{\Real^n},q_2)\in\mathcal{P}(\h),\quad \pr_\h\circ
R(q_1,q_2)=B\in\h.$$ At the same time, the element~$B\in\h$ can be
choosen in an arbitrary way; in order to see that it is enough to
consider the following tensor~$R\in\mathcal{R}(\g^\h)$:
$$R(q_1,q_2)=(B,0,0,0),\quad R(X,Y)=(0,0,0,2\eta(BX,Y)),\quad R(p_1,\cdot)=R(p_2,\cdot)=0,$$
$$R(X,q_1)=(0,0,BX,0),\quad R(X,q_2)=(0,-BX,0,0),\quad
X,Y\in\Real^n.$$ The belonging~$R\in\mathcal{R}(\g^\h)$ can be
checked directly. For any pseudo-Riemannian manifold $(M,g)$ with
the holonomy algebra~$\hol\subset\so(T_xM,g_x)$ at a point~$x\in
M$ we have
$$\nabla_{Z_\alpha}\cdots\nabla_{Z_1}R_x\in\mathcal{R}(\hol),\quad \alpha\geq 0,\quad Z_1,...,Z_\alpha\in T_xM.$$
 In the construction of the metric from the last section we used this property as well as the
 just described algebraic curvature tensors.


\end{document}